\documentclass[12pt]{article}
\usepackage{amssymb}
\usepackage{latexsym,bm}
\usepackage{graphicx}
\usepackage{amsmath}

\newcommand{\bea}{\begin{eqnarray*}}
\newcommand{\eea}{\end{eqnarray*}}
\newcommand{\be}{\begin{equation}}
\newcommand{\ee}{\end{equation}}
\newcommand{\ben}{\begin{eqnarray*}}
\newcommand{\een}{\end{eqnarray*}}

\date{}
\begin{document}
\title{On vertex types of graphs\footnote{E-mail addresses:
{\tt 235711gm@sina.com}(P.Qiao),
{\tt zhan@math.ecnu.edu.cn}(X.Zhan). This research  was supported by the Shanghai SF grant 15ZR1411500 and the NSFC grant 11671148.}}
\author{Pu Qiao\thanks{Corresponding author.}, Xingzhi Zhan\\
{\small Department of Mathematics, East China Normal University, Shanghai 200241, China}
 } \maketitle
\begin{abstract}
 The vertices of a graph are classified into seven types by J.T. Hedetniemi, S.M. Hedetniemi, S.T. Hedetniemi and T.M. Lewis and they ask the following
 questions: 1) What is the smallest order $n$ of a graph having $n-2$ very typical vertices or $n-2$ typical vertices? 2) What is the smallest order of
 a pantypical graph? We answer these two questions in this paper.
\end{abstract}


\section{Introduction}

We consider finite simple graphs. For a vertex $v$ in a graph, we denote by $d(v)$ and $N(v)$ the degree of $v$ and the neighborhood of $v$ respectively
throughout the paper. Motivated by the notions of strong and weak vertices in [3] and [2], J.T. Hedetniemi, S.M. Hedetniemi, S.T. Hedetniemi and T.M. Lewis [1]
classified the vertices of a graph into the following seven types.

{\bf Definition.} A vertex $u$ in a simple graph is said to be
\newline\noindent 1. {\it very strong} if $d(u)\ge 2$ and for every vertex $v\in N(u),$ $d(u)>d(v);$
\newline\noindent 2. {\it strong} if $d(u)\ge 2$ and for every vertex $v\in N(u),$ $d(u)\ge d(v),$ at least one neighbor $x\in N(u)$ has $d(x)<d(u)$
and at least one neighbor $y\in N(u)$ has $d(y)=d(u);$
\newline\noindent 3. {\it regular} if $d(u)\ge 0$ and for every vertex $v\in N(u),$ $d(u)=d(v);$
\newline\noindent 4. {\it very typical} if $d(u)\ge 2$ and for every vertex $v\in N(u),$ $d(u)\ne d(v),$ at least one neighbor $x\in N(u)$ has $d(x)>d(u)$
and at least one neighbor $y\in N(u)$ has $d(y)< d(u);$
\newline\noindent 5. {\it typical} if $d(u)\ge 3$ and there are three distinct vertices $x,y,z\in N(u)$ satisfying $d(x)<d(u)=d(y)<d(z);$
\newline\noindent 6. {\it weak} if $d(u)\ge 2$ and for every vertex $v\in N(u),$ $d(u)\le d(v),$ at least one neighbor $x\in N(u)$ has $d(x)>d(u)$
and at least one neighbor $y\in N(u)$ has $d(y)=d(u);$
\newline\noindent 7. {\it very weak} if $d(u)\ge 1$ and for every vertex $v\in N(u),$ $d(u)< d(v).$

If a graph $G$ has vertices of all seven types, then $G$ is said to be {\it pantypical.}

By the definition above, isolated vertices are regular. Now every simple graph $G$ corresponds to a $7$-tuple $\Gamma(G)=(n_1,n_2,n_3,n_4,n_5,n_6,n_7)$
where $n_1, n_2, n_3, n_4, n_5, n_6, n_7$ are the numbers of very strong, strong, regular, very typical, typical, weak and very weak vertices of $G$ respectively.
We call $\Gamma (G)$ the {\it vertex type} of $G.$  Clearly, two isomorphic graphs must have the same vertex type. Thus the concept of vertex type provides a new
necessary condition for isomorphism when degree sequences cannot distinguish graphs. For example, the two graphs in Figure 1 have the same degree sequence $4,4,4,3,3,2.$
Since the graph (a) has one very weak vertex while the graph (b) has three very weak vertices, the two graphs have different vertex types and hence they are not isomorphic.
\vskip 3mm
\par
 \centerline{\includegraphics[width=4.5in]{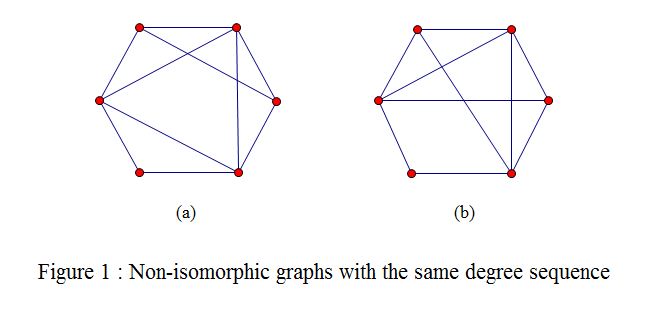}}
\par

The following two questions are asked in [1].
\newline\noindent (1) What is the smallest order $n$ of a graph having $n-2$ very typical vertices or $n-2$ typical vertices?
\newline\noindent (2) What is the smallest order of a pantypical graph?

The purpose of this paper is to answer these two questions. Concerning question (1), since by definition a vertex of maximum degree
or of minimum degree is neither very typical nor typical, a graph of order $n$ can have at most $n-2$ very typical vertices and at
most $n-2$ typical vertices. For every order $n\ge 17$ of the form $n=k^2+1$ with $k$ an integer, a graph of order $n$ with $n-2$
very typical vertices is constructed in [1, proof of Theorem 4], and for every order $n\ge 20$ of the form $n=k(k+1)$ with $k$ an integer, a graph of order $n$ with $n-2$ typical vertices is constructed in [1, proof of Theorem 5].
Concerning question (2), a pantypical graph of order $9$ and size
$21$ is given in [1].

\section{Main Results}

The main results are as follows. All the graphs are simple.

{\bf Theorem 1} {\it Let $f(n)$ and $g(n)$ be the maximum number of very typical vertices and the maximum number of typical vertices
in a graph of order $n.$ Then
$$
f(n)=\begin{cases} 0\quad\quad\,\,\,\, {\rm if}\,\,\, n\le 4;\\
            n-4\quad{\rm if}\,\,\, 5\le n\le 6;\\
            n-3\quad{\rm if}\,\,\, 7\le n\le 9;\\
            n-2\quad{\rm if}\,\,\, n\ge 10
    \end{cases}
$$
and
$$
g(n)=\begin{cases} 0\quad\quad\,\,\,\,{\rm if}\,\,\,n\le 4;\\
                   n-3\quad{\rm if}\,\,\, 5\le n\le 8;\\
                   n-2\quad{\rm if}\,\,\, n\ge 9.
     \end{cases}
$$
}

{\bf Corollary 2} {\it The smallest order $n$ of a graph having $n-2$ very typical vertices is $10$ and the smallest order $n$ of a graph
having $n-2$ typical vertices is $9.$}

{\bf Theorem 3} {\it There exists a pantypical graph of order $n$ if and only if  $n\ge 9.$}

In the following proofs we abbreviate very strong, strong, regular, very typical, typical, weak and very weak as VS, S, R, VT, T, W and VW
respectively. For two vertices $u$ and $v,$ we use the symbol $u\leftrightarrow v$ to mean that $u$ and $v$ are adjacent and use $u\nleftrightarrow v$
to mean that $u$ and $v$ are non-adjacent. The symbol $\Rightarrow$ means ``implies", and $\phi$ denotes the empty set. For two subsets of vertices $P$ and $Q$
in a graph $G$, the symbol $[P,Q]$ denotes the set of those edges with one end vertex in  $P$ and the other end vertex in $Q,$ and $G[P]$ denotes the subgraph
of $G$ induced by $P.$ For a vertex $v,$ $N[v]$  denotes the closed neighborhood of $v;$ i.e., $N[v]=\{v\}\cup N(v).$ Finally $V(G)$ denotes the vertex set
of a graph $G.$

{\bf Proof of Theorem 1.} We first consider $f(n).$ The case $n\le 4$ is trivial. Now suppose $n\ge 5.$ In Figure 2 we give graphs of orders $n=5,6$ with
$n-4$ very typical vertices.
\vskip 3mm
\par
 \centerline{\includegraphics[width=2.8in]{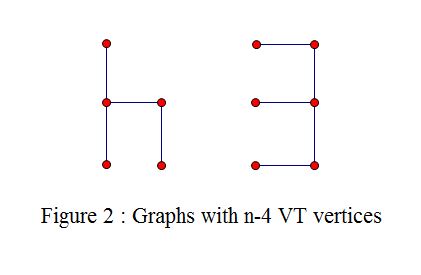}}
\par
In Figure 3 we give graphs of orders $n=7,8,9$ with $n-3$ very typical vertices.
\vskip 3mm
\par
 \centerline{\includegraphics[width=5.8in]{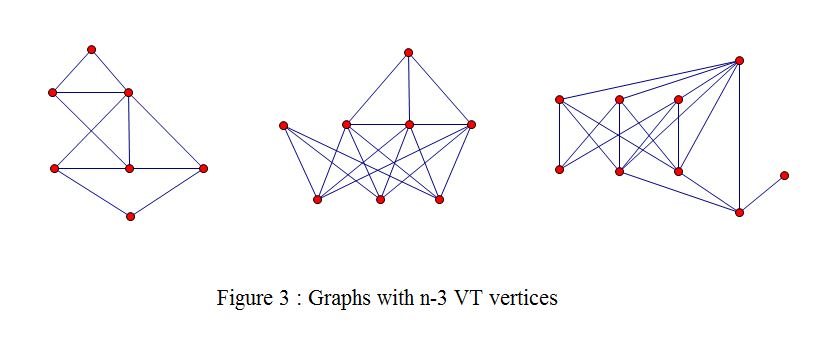}}
\par
In Figure 4 we give graphs of orders $n=10,11$ with $n-2$ very typical vertices.
\vskip 3mm
\par
 \centerline{\includegraphics[width=4.2in]{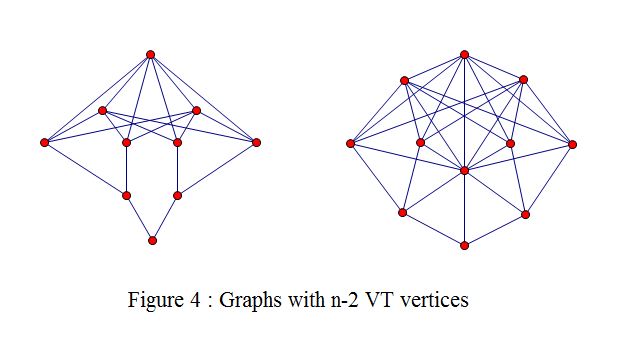}}
\par
If $n=2t\ge 12,$ let $A$ be the complete $4$-partite graph with partite sets of sizes $1,2,t-3$ and $t-1$ respectively. Let $B$ be the graph obtained from $A$
by adding one additional vertex that is adjacent to every vertex in the partite set of size $t-1.$ Then $B$ is a graph of order $n$ with $n-2$ very typical vertices.
If $n=2t+1\ge 13,$ let $C$ be the complete $4$-partite graph with partite sets of sizes $1,2,t-3$ and $t$ respectively. Let $D$ be the graph obtained from $C$ by adding
one additional vertex that is adjacent to every vertex in the partite set of size $t.$ Then $D$ is a graph of order $n$ with $n-2$ very typical vertices.

Since a vertex of maximum or minimum degree cannot be a very typical vertex, a graph of order $n$ has at most $n-2$ very typical vertices. Thus, the above constructions
show that for $5\le n\le 6,$ $f(n)\ge n-4,$ for $7\le n\le 9,$ $f(n)\ge n-3$ and
for $n\ge 10,$ $f(n)=n-2.$ Next we show that a graph of order $9$ cannot have $n-2=7$ very typical vertices. The proofs that a graph of orders $n=7,8$ cannot have
$n-2$ very typical vertices and that a graph of orders $n=5,6$ cannot have $n-3$ very typical vertices are similar to the case $n=9,$ but are easier, so we omit them and assume
that the results for these lower orders have been proved.

To the contrary, suppose that there is a graph $G$ of order $9$ with $7$ VT vertices. Let $V(G)=\{v_1,v_2,\ldots,v_9\}$ with $d(v_1)> d(v_2)\ge\cdots\ge d(v_8)>d(v_9).$
Then $v_2,\ldots,v_8$ are VT. Clearly $d(v_9)=0$ is impossible. If $d(v_9)=1,$ then $G-v_9$ is a graph of order $8$ with $6$ VT vertices, a contradiction. From
now on, we assume that $d(v_9)\ge 2.$ Denote by $S_i$ the set of vertices in $V(G)\setminus \{v_1,v_9\}$ with degree $i$ for $i=3,4,5,6,7.$
Then each $S_i$ is an independent set.
We have the rough estimate that
$$
|S_7|\le 2,\,\,\, |S_6|\le 3,\,\,\, |S_5|\le 4,\,\,\, |S_4|\le 5,\,\,\, |S_3|\le 2.
$$
Note that the vertices in $S_3$ have $v_9$ as the only lower degree neighbor. If $|S_3|\ge 3,$ then $d(v_9)\ge 3,$ implying that $S_3=\phi,$ a contradiction.
Thus we have $|S_3|\le 2.$

Now we show that $S_7=\phi.$ If $|S_7|=2,$ then $d(v_1)=8,$ $d(v_9)\ge 3$ and $S_7=\{v_2,v_3\}.$ But then $G-\{v_1,v_2\}$ is a graph of order $7$ with $5$ VT vertices,
a contradiction. If $|S_7|=1,$ then $d(v_1)=8$ and $S_7=\{v_2\}.$ Let $w$ be the vertex such that $w\nleftrightarrow v_2.$ If $w$ has a higher degree neighbor other than $v_1,$
then $G-v_1$ is a graph of order $8$ with $6$ VT vertices, a contradiction. If $v_1$ is the only higher degree neighbor of $w,$ then in $G-v_1$ we add the edge
$wv_2$ to obtain a graph of order $8$ with $6$ VT vertices, a contradiction.

If $|S_6|=3,$ then $d(v_9)\ge 4$ $\Rightarrow$ $S_3=S_4=\phi$ $\Rightarrow$ $|S_5|=4$ $\Rightarrow$ $d(v_9)\ge 7,$ a contradiction. Thus $|S_6|\le 2.$ If $|S_5|=4,$ then
$d(v_9)\ge 4$ $\Rightarrow$ $S_3=S_4=\phi$ $\Rightarrow$ $|S_6|=3,$ impossible. If $|S_4|=5,$ then $d(v_9)\ge 5$ $\Rightarrow$ $S_3=S_4=S_5=\phi,$ a contradiction. Thus
we have the following sharper estimate that
$$
S_7=\phi,\,\,\, |S_6|\le 2,\,\,\, |S_5|\le 3,\,\,\, |S_4|\le 4,\,\,\, |S_3|\le 2.
$$

Now we show that $S_6=\phi.$ First suppose $|S_6|=2.$ Then $S_6=\{v_2,v_3\}.$ We distinguish three cases.

Case 1. $v_2\nleftrightarrow v_9$ and $v_3\nleftrightarrow v_9.$ In this case $v_2$ and $v_3$ are adjacent to each vertex in $S_5\cup S_4\cup S_3$ and $S_6\neq\phi$
$\Rightarrow$ $d(v_1)\ge 7.$ But then either at least one vertex in $S_5\cup S_4$ has no lower degree neighbor or $d(v_1)\le 5,$ a contradiction.

Case 2. $v_2\leftrightarrow v_9$ and $v_3\leftrightarrow v_9.$  We have $S_3=\phi$ $\Rightarrow$ $|S_4\cup S_5|=5.$ But $|S_5|\le 3$ $\Rightarrow$ $|S_4|\ge 2$
$\Rightarrow$ $d(v_9)\ge 4$ $\Rightarrow$ $S_4=\phi,$ a contradiction.

Case 3. One of $v_2$ and $v_3$, say $v_2\leftrightarrow v_9$ and $v_3\nleftrightarrow v_9.$ We have $|S_3|\le 1.$ Note that $v_3$ is adjacent to each vertex in
$S_3\cup S_4\cup S_5.$ If $|S_3|=1,$ then $|S_4|=1.$ To see this, note that $|S_5|\le 3$ $\Rightarrow$ $|S_4|\ge 1,$ but $|S_4|>1$ would imply that at least one vertex
in $S_4$ has no lower degree neighbor.
Consequently $|S_5|=3,$ $S_3=\{v_8\}$ and $S_4=\{v_7\}.$ Now since the three vertices in
$S_5$ can have $v_7$ as the only lower degree neighbor, we have $d(v_7)\ge 5,$ a contradiction. If $|S_3|=0,$ then $|S_5|\le 3$ $\Rightarrow$ $|S_4|\ge 2$
$\Rightarrow$ $d(v_9)\le 3$ $\Rightarrow$ $|S_4|\le 2.$ Hence $|S_4|=2$ and $d(v_9)=3$ $\Rightarrow$ $S_4=\{v_8,v_7\}.$ Now the three vertices in $S_5$
can have either $v_8$ or $v_7$ as a lower degree neighbor and $v_2$ is adjacent to at least one of $v_8$ and $v_7.$ Hence $N(v_1)\subseteq S_6\cup S_5$ $\Rightarrow$
$d(v_1)\le 5,$ a contradiction.

Next suppose $|S_6|=1.$ Then $S_6=\{v_2\}.$ We distinguish two cases.

Case 1. $v_2\leftrightarrow v_9.$ Then we have $|S_3|\le 1.$

Subcase 1. $|S_3|=1.$ This implies that $S_3=\{v_8\}$ and $|S_4|\le 2.$ But $|S_5|\le 3$ and $|S_5|+|S_4|=5.$ Hence $|S_4|=2.$ We deduce that $v_9,v_8,v_1\notin N(v_1)$
$\Rightarrow$ $d(v_1)\le 6,$ contradicting $d(v_1)\ge 7.$

Subcase 2. $|S_3|=0.$ Now $|S_4|+|S_5|=6.$  $S_3=\phi$ and $S_4\neq\phi$ $\Rightarrow$ $d(v_9)\le 3$ $\Rightarrow$ $|S_4|\le 2.$ Hence $|S_5|\ge 4,$ contradicting $|S_5|\le 3.$

Case 2. $v_2\nleftrightarrow v_9.$ Since $|S_3|\le 2,$ we distinguish the following three subcases.

Subcase 1. $|S_3|=0.$ We have $|S_4|+|S_5|=6$ and $|S_5|\le 3$ $\Rightarrow$ $|S_4|\ge 3$ $\Rightarrow$ $d(v_9)\le 3$ $\Rightarrow$ $|S_4|=3$ $\Rightarrow$ $|S_5|=3$
and $d(v_9)=3.$ Let $S_5=\{a_1,a_2,a_3\}$ and $S_4=\{b_1,b_2,b_3\}.$ Note that $N(v_9)=S_4.$ Since each $a_i\nleftrightarrow v_9,$ $S_5$ is independent and $d(a_i)=5,$
we deduce that each $a_i \leftrightarrow$ each $b_j$ $\Rightarrow$ $N(b_j)=\{v_9,a_1,a_2,a_3\}$ for every $j=1,2,3$ $\Rightarrow$ $N(v_2)\subseteq \{v_1\}\cup S_5$
$\Rightarrow$ $d(v_2)\le 4,$ contradicting $d(v_2)=6.$

Subcase 2. $|S_3|=1.$ $S_3\neq \phi$ $\Rightarrow$ $d(v_9)\le 2.$ We have $|S_4|\le 3,$ since otherwise at least one vertex in $S_4$ would have no lower degree neighbor.
$|S_4|+|S_5|=5$ and $|S_5|\le 3$ $\Rightarrow$ $|S_4|\ge 2.$ Thus there are two possibilities: $|S_4|=2$ or $3.$

Let $P=\{v_1\}\cup S_6\cup S_5$ and $Q=S_4\cup S_3\cup \{v_9\}.$ Note that $d(v_1)\ge 7$ and $d(v_9)=2.$

If $|S_4|=2,$ then $|S_5|=3$ and
$$
\sum_{v\in P}d(v)-\sum_{v\in Q}d(v)\ge 7+6+3\times 5-2\times 4-3-2=15.\eqno (1)
$$
Since every edge in  $[P,Q]$ contributes the same degree $1$ to both $\sum_{v\in P}d(v)$ and $\sum_{v\in Q}d(v),$ to calculate their difference it suffices to
consider the degrees from those edges inside $G[P]$ or $G[Q].$ There are at least three edges inside $G[Q].$ Thus
$$
\sum_{v\in P}d(v)-\sum_{v\in Q}d(v)\le 5\times 4-3\times 2=14,
$$
contradicting (1).

If $|S_4|=3,$ then $|S_5|=2$ and
$$
\sum_{v\in P}d(v)-\sum_{v\in Q}d(v)\ge 7+6+2\times 5-3\times 4-3-2=6.\eqno (2)
$$
On the other hand, since there are at least four edges inside $G[Q]$ we have
$$
\sum_{v\in P}d(v)-\sum_{v\in Q}d(v)\le 4\times 3-4\times 2=4,
$$
contradicting (2).

Subcase 3. $|S_3|=2.$ $S_3\neq \phi$ $\Rightarrow$ $d(v_9)=2.$ $|S_4|+|S_5|=4$ and $|S_5|\le 3$ $\Rightarrow$ $|S_4|\ge 1.$ Since $d(v_1)\ge 7$ and
$v_1\nleftrightarrow v_9,$ each vertex in $S_3$ is adjacent to $v_1.$ Consequently $|S_4|\le 2.$ Using the same method as in the above subcase 2 to the
two possible cases $|S_4|=1$ and $|S_4|=2$ we obtain contradictions too.

So far we have proved that $S_7=S_6=\phi.$ Next according to $|S_3|\le 2$ we distinguish three cases.

Case 1. $|S_3|=2.$ We have $S_3=\{v_8,v_7\},d(v_9)=2,$ $|S_4|+|S_5|=5$ and $|S_5|\le 3$ $\Rightarrow$ $|S_4|\ge 2.$ $|S_4|\le 4$ $\Rightarrow$
$S_5\neq \phi$ $\Rightarrow$ $d(v_1)\ge 6$
$\Rightarrow$ $v_1$ is adjacent to at least one of $v_8$ and $v_7$ $\Rightarrow$ $2\le |S_4|\le 3.$ If $|S_4|=2,$ then $|S_5|=3$ and at least one vertex in $S_5$
has degree $\le 3,$ a contradiction. If $|S_4|=3,$ then $|S_5|=2$ and the two vertices in $S_5$ have degree $\le 4,$ a contradiction.

Case 2. $|S_3|=1.$ We have $S_3=\{v_8\}.$ The vertices in $S_4$ can only have $v_8$ or $v_9$ as their lower degree neighbors. Since $d(v_9)=2,$ $|S_4|\le 3.$
But $|S_4|+|S_5|=6$ and $|S_5|\le 3.$ Hence $|S_4|=|S_5|=3.$ Then the vertices in $S_5$ have degree $\le 4,$ a contradiction.

Case 3. $|S_3|=0.$ Now $|S_4|+|S_5|=7,$ $|S_4|\le 4,$ $|S_5|\le 3$ $\Rightarrow$ $|S_4|=4.$ Since the vertices in $S_4$ can only have $v_9$ as their lower degree
neighbor, they must be adjacent to $v_9.$ Hence $d(v_9)\ge 4.$ But then the vertices in $S_4$ have no lower degree neighbor, a contradiction. This completes the
proof of the result on $f(n).$

Now we consider $g(n).$ The case $n\le 4$ is trivial and we assume $n\ge 5.$ In Figure 5 we give graphs of orders $n=5,6,7,8$ with $n-3$ T vertices.
\vskip 3mm
\par
 \centerline{\includegraphics[width=5.8in]{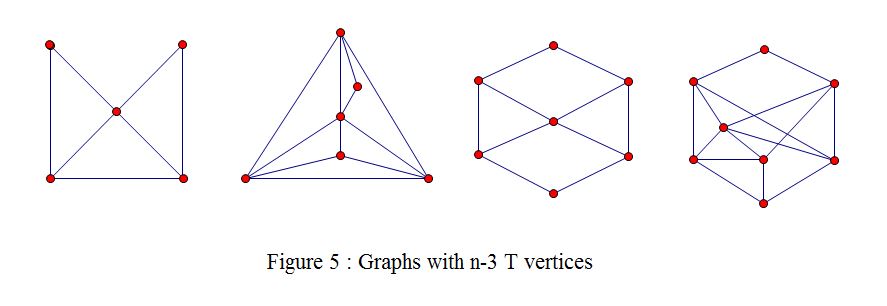}}
\par
For each $n\ge 9,$ we construct a graph of order $n$ with $n-2$ T vertices. Recall that the join of graphs $G_1,G_2,\ldots,G_k,$ denoted
$G_1\vee G_2\vee\cdots \vee G_k,$ is the graph obtained from the disjoint union $G_1+G_2+\cdots +G_k$ by adding the edges $xy$ with $x\in V(G_i)$ and
$y\in V(G_j)$ for all $i\neq j.$

In the following constructions we let $K_1$ be the graph consisting of one vertex, let $C_t$ be the cycle of order $t,$ let $M_{2h}$ be a graph of order $2h$
and size $h$ whose edges form a perfect matching, and let $T_s$ be a cubic graph of order $s.$

If $n=4k+1\ge 9,$ to the graph $K_1\vee C_{2k-1}\vee M_{2k}$ add one additional vertex that is adjacent to each vertex in $V(K_1)\cup V(M_{2k}).$ Then we obtain
a graph of order $n$ with $n-2$ T vertices.

If $n=4k+2\ge 10,$ to the graph $K_1\vee T_{2k}\vee M_{2k}$ add one additional vertex that is adjacent to each vertex in $V(K_1)\cup V(M_{2k}).$ Then we obtain
a graph of order $n$ with $n-2$ T vertices.

If $n=4k+3\ge 11,$ to the graph $K_1\vee T_{2k}\vee C_{2k+1}$ add one additional vertex that is adjacent to each vertex in $V(K_1)\cup V(C_{2k+1}).$ Then we obtain
a graph of order $n$ with $n-2$ T vertices.

If $n=4k\ge 12,$ to the graph $K_1\vee M_{2k-2}\vee M_{2k}$ add one additional vertex that is adjacent to each vertex in $V(M_{2k}).$ Then we obtain
a graph of order $n$ with $n-2$ T vertices.

Since a vertex of maximum or minimum degree cannot be a T vertex, a graph of order $n$ can have at most $n-2$ T vertices. Thus the above constructions show that
 for $5\le n\le 8,$ $g(n)\ge n-3$ and for $n\ge 9,$ $g(n)=n-2.$ It remains to prove that for $5\le n\le 8,$ $g(n)\le n-3.$ The proofs for the cases $5\le n\le 7$
 are similar to that for the case $n=8$ but easier, so we omit them.

 Let $G$ be a graph of order $8$ with $V(G)=\{v_1,v_2,\ldots,v_8\}$ and $d(v_1)\ge d(v_2)\ge\cdots\ge d(v_8).$ We will show that $G$ has at most $5$ T vertices.
 To the contrary, suppose $G$ has $6$ T vertices. Then $v_2,\ldots, v_7$ are T vertices and
 $$
 d(v_1)>d(v_2)=d(v_3)\ge d(v_4)\ge d(v_5)\ge d(v_6)=d(v_7)>d(v_8).
 $$
Renaming the vertices $v_3,v_4,v_5$ and $v_6$ if necessary, we may assume that $v_2\leftrightarrow v_3$ and  $v_6\leftrightarrow v_7.$
Then $v_1v_2v_3$ and $v_6v_7v_8$ are triangles. Denote $a=d(v_2)=d(v_3)$ and $b=d(v_6)=d(v_7).$ Clearly $a>b.$ We have
$b\ge 3$ $\Rightarrow$ $a\ge 4$ $\Rightarrow$ $d(v_1)\ge 5.$ Also, $d(v_1)\le 7$ $\Rightarrow$ $a\le 6$ $\Rightarrow$ $b\le 5$ $\Rightarrow$ $d(v_8)\le 4.$
We assert that $d(v_4)>b.$ Otherwise $d(v_4)=d(v_5)=b$ $\Rightarrow$ $v_4,v_5,v_6,v_7\in N(v_8)$ $\Rightarrow$ $d(v_8)\ge 4$ $\Rightarrow$ $b\ge 5$
$\Rightarrow$ $a\ge 6$ $\Rightarrow$ $d(v_1)=7$ $\Rightarrow$ $v_1\leftrightarrow v_8$ $\Rightarrow$ $d(v_8)\ge 5,$ contradicting $d(v_8)\le 4.$
We further assert that $d(v_5)<a.$ Otherwise $d(v_4)=d(v_5)=a.$
Let $P=\{v_1,v_2,v_3,v_4,v_5\}$ and $Q=\{v_6,v_7,v_8\}.$ $d(v_1)\ge 5$ $\Rightarrow$ $|[v_1,Q]|\ge 1.$ Every vertex in $\{v_2,v_3,v_4,v_5\}$ has at least
one lower degree neighbor in $Q.$ Hence $|[P,Q]|\ge 1+4=5,$ which, together with the fact that $v_6v_7v_8$ is a cycle, implies that
$\sum_{v\in Q}d(v)\ge 5+3\times 2=11.$ Since $b>d(v_8),$ $b\ge 4$ $\Rightarrow$ $a\ge 5$ $\Rightarrow$ $d(v_1)\ge 6$ $\Rightarrow$ $|[P,Q]|\ge 2+4=6$
$\Rightarrow$ $\sum_{v\in Q}d(v)\ge 6+3\times 2=12$ $\Rightarrow$ $b\ge 5$ $\Rightarrow$ $a\ge 6$ $\Rightarrow$ $a=6$ and $d(v_1)=7$ $\Rightarrow$
$|[v_1,Q]|\ge 3$ and $|[v_i,Q]|\ge 2$ for $i=2,3,4,5$ $\Rightarrow$ $|[P,Q]|\ge 3+4\times 2=11$ $\Rightarrow$ $\sum_{v\in Q}d(v)\ge 11+3\times 2=17$
$\Rightarrow$ $b\ge 6,$ contradicting $b\le 5.$

Next we distinguish two cases.

Case 1. $d(v_5)>b.$ Then $a>d(v_4)=d(v_5)>b.$ Let $P=\{v_1,v_2,v_3\}$ and $Q=\{v_6,v_7,v_8\}.$ Since $d(v_1)\ge d(v_8)+4,$ $d(v_2)\ge d(v_6)+2$
and $d(v_3)\ge d(v_7)+2,$ we have $\sum_{v\in P}d(v)-\sum_{v\in Q}d(v)\ge 8.$ But this is impossible. Since $v_1v_2v_3$ and $v_6v_7v_8$ are both triangles,
more degrees in  $\sum_{v\in P}d(v)$ than in $\sum_{v\in Q}d(v)$ can only come from edges incident to $v_4$ and $v_5$ and hence
$\sum_{v\in P}d(v)-\sum_{v\in Q}d(v)\le 3\times 2=6.$

Case 2.  $d(v_5)=b.$ Then $d(v_4)=a$ since $d(v_4)>b,$ $\{v_2,v_3,v_4\}\subseteq N(v_1)$ and $\{v_5,v_6,v_7\}\subseteq N(v_8).$ Let $P=\{v_1,v_2,v_3,v_4\}$
and $Q=\{v_5,v_6,v_7,v_8\}.$  Since $d(v_1)\ge d(v_8)+3,$ $d(v_2)\ge d(v_5)+1,$ $d(v_3)\ge d(v_6)+1$ and $d(v_4)\ge d(v_7)+1,$ we have
$\sum_{v\in P}d(v)-\sum_{v\in Q}d(v)\ge 6.$ On the other hand, since each of the induced subgraphs $G[P]$ and $G[Q]$ has size either $5$ or $6,$
$\sum_{v\in P}d(v)-\sum_{v\in Q}d(v)\le 2,$ a contradiction. This completes the proof.$\Box$

Corollary 2 follows from Theorem 1 immediately. We will repeatedly use the following lemma which follows from the definition.

{\bf Lemma 4} Let $s^{\prime},\,s,\,r,\,t^{\prime},\,w$ and $w^{\prime}$ be a $VS,\,S,\,R,\,VT,\,W$ and $VW$ vertex respectively in a graph. Then
$$
s^{\prime}\nleftrightarrow s,\,\,\,w^{\prime}\nleftrightarrow w,\,\,\,r\nleftrightarrow t^{\prime},\,\,\,r\nleftrightarrow s^{\prime},\,\,\,r\nleftrightarrow w^{\prime}.
$$

{\bf Proof of Theorem 3.} In Figure 6 we give a pantypical graph of order $9$ and size $11.$
\vskip 3mm
\par
 \centerline{\includegraphics[width=3.8in]{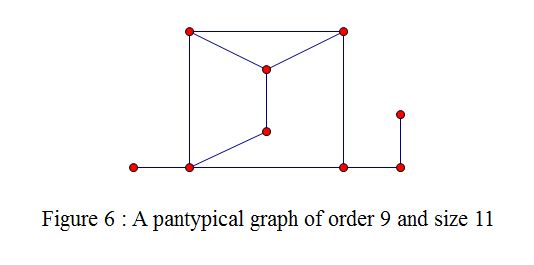}}
\par
For $n>9,$ we attach a path of order $n-8$ to the vertex of maximum degree in Figure 6 to obtain a pantypical graph of order $n.$

Conversely we need to prove that there is no pantypical graph of order $\le 8.$ Obviously a pantypical graph must have order at least $7.$
Since the proof for the case $n=7$ is similar but easier, we present only the proof for the case $n=8.$

To the contrary, suppose that there is a pantypical graph $G$ of order $8.$
Since a VT, T, W or VW vertex has a higher degree neighbor, it cannot have the maximum degree. Hence $G$ has at most $4$ vertices of the maximum degree.
Let $\Delta$ be the maximum degree of $G.$ Then $4\le\Delta\le 7.$ We
distinguish four cases according to the values of $\Delta.$

$$
\Delta =7
$$

By considering a vertex of degree $7$ and a regular vertex we conclude that $G=K_8$ which is not pantypical, a contradiction.

$$
\Delta =6
$$

Let $v$ be a vertex of $G$ with degree $6.$ Then $v$ is VS or S. If $v$ is VS, then $v$ has a neighbor which is S or R, contradicting Lemma 4. If $v$ is S,
then the only vertex $\notin N[v]$ must be VS and consequently a regular vertex is adjacent to $v.$ It follows that $G$ has at least $7$ vertices
of the maximum degree $6,$ which is impossible.

$$
\Delta =5
$$

Let $v$ be a vertex of $G$ with degree $5,$ and let $r$ be a regular vertex. Then $r\notin N(v),$ since otherwise $G$ would have at least $6$ vertices
with the maximum degree $5,$ which is impossible. Note that $v$ must be either S or VS.

Case 1. $v$ is S. Let $s^{\prime}$ be a VS vertex. Then $s^{\prime}\nleftrightarrow v.$ By definition, $v$ has a neighbor $x$ with $d(x)=5.$ Then $x$ is S. Let
$N(v)=\{x,y,z,t,t^{\prime}\}$ where $t$ is T and $t^{\prime}$ is VT. Clearly $N(x)=\{v,y,z,t,t^{\prime}\}.$ Since $d(t)\ge 3$ and $t$ already has two neighbors $v$ and $x$
of the maximum degree $5,$ $d(t)=4.$ It is easy to see that $t\nleftrightarrow r$ and $t\nleftrightarrow s^{\prime}.$ Now $v,t,x,r\notin N(s^{\prime})$ and hence
$y,z,t^{\prime}$ are the only possible neighbors of $s^{\prime}$ $\Rightarrow$ $d(s^{\prime})\le 3.$
Since $v,x\in N(y)\cap N(z)\cap N(t^{\prime}),$ $s^{\prime}$ cannot have any lower degree neighbors, a contradiction.

Case 2. $v$ is VS. In this case, among the two non-neighbors of $v,$ one is R, denoted $r$ and the other is S, denoted $s.$ Let $N(v)=\{t,t^{\prime},a,b,c\}$
where $t$ is T and $t^{\prime}$ is VT. We have $3\le d(t)\le 4.$ Note that $d(t^{\prime})\ge 3,$ since any lower degree neighbor of $t^{\prime}$ has degree
$\ge 2.$

Subcase 1. $d(t)=3.$
Since $d(t^{\prime})\ge 3$ and $t$ already has a higher degree neighbor $v,$ $t\nleftrightarrow t^{\prime}.$ The lower degree neighbor of $t$ is in $\{a,b,c\},$
say $t\leftrightarrow a$ and $d(a)=2.$ Clearly $t$ cannot be adjacent to both $r$ and $s.$ Thus we have three cases.
(1) $t\leftrightarrow r.$ If $r\leftrightarrow s,$ then $d(r)=d(s)=d(t)=3.$ Since $d(t^{\prime})\ge 3,$ $s$ has at least one lower degree neighbor in $\{b,c\},$
say $s\leftrightarrow b$ and $d(b)=2.$ Then $d(r)=3$ $\Rightarrow$ $r\leftrightarrow c.$ But now $t^{\prime}$ has no lower degree neighbor, a contradiction.
If $r\nleftrightarrow s,$ then $d(r)=3$ $\Rightarrow$ $r\leftrightarrow b$ and $r\leftrightarrow c.$
Now $t^{\prime}$ can only have a lower degree neighbor in $\{b,c\}$ which has degree $3$ $\Rightarrow$ $d(t^{\prime})=4$
$\Rightarrow$ $t^{\prime}\leftrightarrow b,$ $t^{\prime}\leftrightarrow c$ and $t^{\prime}\leftrightarrow s$ $\Rightarrow$ $d(s)=1,$
 a contradiction.
(2) $t\leftrightarrow s.$ Now $d(s)=3$ $\Rightarrow$ $s$ has at least one lower degree neighbor in $\{b,c\},$ say $s\leftrightarrow b$ and $d(b)=2.$
Consequently $N(t^{\prime})\subseteq \{v,c,s\}$ $\Rightarrow$ $d(t^{\prime})\le 3$ $\Rightarrow$ $d(t^{\prime})=3.$ But $d(s)=3$ $\Rightarrow$
$t^{\prime}\nleftrightarrow s.$ Then $d(t^{\prime})\le 2,$ a contradiction.
(3) $t\nleftrightarrow r$ and $t\nleftrightarrow s.$ In this case, the neighbor of $t$ with degree $3$ can only be $b$ or $c,$ say $t\leftrightarrow b$ and $d(b)=3.$
We have $N(t^{\prime})\subseteq \{v,b,c,s\}$ $\Rightarrow$ $d(t^{\prime})=3$ or $4$ $\Rightarrow$ $t^{\prime}\leftrightarrow s$ $\Rightarrow$ $d(s)\ge 4.$ Since
$N(s)\subseteq \{r,b,c,t^{\prime}\}$ we have $d(s)=4$ and $N(s)=\{r,b,c,t^{\prime}\}$ $\Rightarrow$ $d(t^{\prime})=3$ $\Rightarrow$ $t^{\prime}\leftrightarrow c$
$\Rightarrow$ $d(c)\ge 3,$ but now $t^{\prime}$ has no lower degree neighbor, a contradiction.

Subcase 2. $d(t)=4.$ Suppose $t\leftrightarrow r.$ Then $d(r)=4.$ But $r$ is not adjacent to any VW or VT vertex by Lemma 4. Without loss of generality,
 suppose $c$ is VW. Then $N(r)=\{t,a,b,s\},$ and $r,t,a,b,s$ all have degree $4.$ Now $t^{\prime}$ can only have $c$ as its lower degree neighbor $\Rightarrow$
$c\leftrightarrow t^{\prime}$ and $d(c)\ge 2.$ But then $t$ and $s$ cannot both have a lower degree neighbor, a contradiction.

Now suppose $t\nleftrightarrow r.$ (1) Suppose  $t\leftrightarrow t^{\prime}.$ Then $d(t^{\prime})=3$ and $t^{\prime}$ has a lower degree neighbor  in $\{a,b,c\},$
say $c,$  of degree $2.$ The vertex $t$ has a neighbor $w$ in $\{a,b,s\}$ of degree $4,$ which has possible neighbors in $\{t,a,b,r,s\}$ except itself.
Hence $w\leftrightarrow r$ and $d(r)=4.$ But $r$ has only the three possible neighbors $a,b,s,$ a contradiction. (2) Suppose  $t\nleftrightarrow t^{\prime}.$
If $t\nleftrightarrow s,$
then $N(t)=\{v,a,b,c\}$ $\Rightarrow$ At least one of $\{a,b,c\}$ has degree $4$ and any lower degree neighbor of $t^{\prime}$ has degree $\ge 3.$
Hence $d(t^{\prime})=4$ $\Rightarrow$ $t^{\prime}\leftrightarrow s$ $\Rightarrow$ $d(s)=5$ $\Rightarrow$ $N(s)=\{r,a,b,c,t^{\prime}\}$ $\Rightarrow$
$d(r)=5,$ contradicting the fact that $N(r)\subseteq \{s,a,b,c\}$ now.
If $t\leftrightarrow s,$ then $d(s)\ge 4.$ Clearly $r\nleftrightarrow s,$
since otherwise $a,b,c\in N(r)$ and they all have degree $\ge 4,$ implying that $t$ has no lower degree neighbor. If $s\nleftrightarrow t^{\prime},$ then
$N(s)=\{t,a,b,c\},$ implying that $t^{\prime}$ has no lower degree neighbor. Hence $s\leftrightarrow t^{\prime}$ $\Rightarrow$ $d(t^{\prime})\neq d(s)$
$\Rightarrow$ $d(t^{\prime})\le 3,$ since $d(s)\ge 4$ and $d(t^{\prime})\le \Delta -1=4$
$\Rightarrow$
 $d(t^{\prime})=3,$ since $d(t^{\prime})\ge 3.$
 Consequently $t^{\prime}$ has a lower degree
neighbor in $\{a,b,c\},$ say $c,$ of degree $2$ $\Rightarrow$ $s\nleftrightarrow c,$ $t\leftrightarrow a$ and $t\leftrightarrow b$
 $\Rightarrow$ $N(s)=\{t,t^{\prime},a,b\}$ $\Rightarrow$ $r\nleftrightarrow a$
and $r\nleftrightarrow b$ $\Rightarrow$ $d(a)=d(b)=3.$ Finally $G$ has three VW vertices $a,b,c,$ a contradiction.

$$
\Delta =4
$$

Let $v$ be a vertex of $G$ with degree $4.$ Then $v$ is S or VS.

Case 1. $v$ is S. $v$ has a neighbor $x$ with $d(x)=4.$ It is easy to see that an R or VS vertex is not in $N(v)\cup N(x).$ Since $G$ has two S vertices $v$ and $x,$
$G$ has exactly one vertex of each of the remaining six types. Let $t,r,s^{\prime}$ be the T, R, VS vertex of $G$ respectively.

 Subcase 1. $t\in N(v).$ Then $d(t)=3,$ $t\nleftrightarrow x$ and $t\nleftrightarrow s^{\prime}.$ Let $N(v)=\{x,y,z,t\}$ and $V(G)\setminus N[v]=\{a,r,s^{\prime}\}.$
 Note that $N(x)=\{v,a,y,z\}.$
 Now $a$ is the lower degree neighbor of $t$ $\Rightarrow$ $t\leftrightarrow a$ and $d(a)=2.$
 If $t\leftrightarrow r,$ then $d(r)=3$   $\Rightarrow$ $r\leftrightarrow y$ and $r\leftrightarrow z$
 $\Rightarrow$ $d(y)=d(z)=3.$ But now $s^{\prime}$ has no neighbors, a contradiction.  Thus, $t\nleftrightarrow r.$

 As a T vertex, $t$ must have a neighbor of the same degree and a neighbor of lower degree.
We deduce that $t\leftrightarrow a,$ $d(a)=2$ and that $t$ is adjacent to one of $y$ and $z,$ say $z:\,t\leftrightarrow z.$  Hence $d(z)=3.$
But now $y$ is the only possible neighbor of $s^{\prime}$ $\Rightarrow$ $d(s^{\prime})\le 1,$ a contradiction.

Subcase 2. $t\notin N(v).$ We have $V(G)\setminus N[v]=\{t,r,s^{\prime}\}.$ Let $N(v)=\{x,y,z,b\}.$ Suppose $x\leftrightarrow t.$ Then $d(t)=3.$ If $t\leftrightarrow r,$
then $d(r)=3.$ Since among the three vertices $y,z,b,$ one is VW and one is VT, neither of which is adjacent to $r.$ Also $r\nleftrightarrow s^{\prime}.$ Hence $d(r)\le 2,$
a contradiction. This shows $t\nleftrightarrow r.$ Without loss of generality, suppose besides $v$ and $t,$ the other two neighbors of $x$ are $y$ and $z.$
Note that
$t\nleftrightarrow s^{\prime}$ since $d(t)=3$ and $t$ already has a higher degree neighbor $x.$
Then $t\leftrightarrow b$ and $t$ is adjacent to exactly one of $y$ and $z,$ say $z.$ Now $N(t)=\{x,z,b\}$ $\Rightarrow$ $d(z)=3$ and $d(b)=2.$
 Then $s^{\prime}$ has $y$ as the only possible neighbor $\Rightarrow$ $d(s^{\prime})\le 1,$
a contradiction. Hence $x\nleftrightarrow t.$ Now $N(x)=\{v,y,z,b\}.$ One of $y,z$ and $b$ is VT with degree $\le 3,$ but it has no lower degree neighbor now, a contradiction.

Case 2. $v$ is VS. By Lemma 4, an R or S vertex is not adjacent to $v.$ We distinguish two subcases according as whether $v$ has a neighbor which is VT.

Subcase 1. $v$ has a neighbor $t^{\prime}$ which is VT. Let $t$ be a T vertex. Clearly $d(t)=3.$ First suppose $d(t^{\prime})=2.$ Then $t^{\prime}$ has a neighbor $w^{\prime}$ of degree 1, and
$V(G)\setminus N[v]=\{w^{\prime},s,r\}$ where $s$ is strong and $r$ is regular. Let $N(v)=\{t^{\prime},a,b,t\}.$  Now both $t$ and $s$ can only have lower degree neighbors
in $\{a,b\}.$ Hence $d(s)\ge 3.$ Consequently $s\leftrightarrow t$ and $s\leftrightarrow r.$ But  $r$ can have $s$ as the only neighbor $\Rightarrow$ $d(r)=1,$ which is a
contradiction since $r$ is regular, $r\leftrightarrow s$ and $d(s)\ge 3.$

Next suppose $d(t^{\prime})=3.$ Since $d(t)=3,$ $t^{\prime}\nleftrightarrow t.$ If $t\nleftrightarrow v,$ let $N(v)=\{t^{\prime},a,b,c\}.$ Both $t^{\prime}$ and $t$ have a lower degree neighbor in
$\{a,b,c\}.$ Without  loss of generality, suppose $t^{\prime}\leftrightarrow a$ and $t\leftrightarrow b.$ We have $d(a)=d(b)=2.$ Now $s$ is the only possible higher degree
neighbor of $t.$ Hence $s\leftrightarrow t$ and $d(s)=4$ $\Rightarrow$ $s\leftrightarrow t^{\prime},$ $s\leftrightarrow c$ and $s\leftrightarrow r$ $\Rightarrow$
$d(r)=4.$ On the other hand, $N(r)\subseteq \{s,t,c\}$ $\Rightarrow$ $d(r)\le 3,$ a contradiction.

If $t\leftrightarrow v,$ let $N(v)=\{t^{\prime},t,x,y\}$ and $V(G)\setminus N[v]=\{z,r,s\}$ where $r$ is R and $s$ is S. We first assert that $t^{\prime}\nleftrightarrow s.$
Otherwise $t^{\prime}\leftrightarrow s$ $\Rightarrow$ $d(s)=4,$ $s\nleftrightarrow t,$ since $d(t)=3$ and $t$ already has a higher degree neighbor $v.$ By definition, $s$
has a neighbor with degree $4,$ which can only be $r$ or $z.$ If $s\leftrightarrow r,$ then $d(r)=4.$ But $r$ has no so many neighbors of the same degree. If
$s\leftrightarrow z$ and $d(z)=4,$ then $z\nleftrightarrow t$  $\Rightarrow$  $z\leftrightarrow x,$ $z\leftrightarrow y,$ $z\leftrightarrow t^{\prime}.$ But then
$t^{\prime}$ has no lower degree neighbor.

If $t^{\prime}\leftrightarrow z,$ then $t^{\prime}\leftrightarrow$ one of $x$ and $y,$ say $t^{\prime}\leftrightarrow x$ $\Rightarrow$ $d(x)=2.$ $y$ and $z$ are the only possible lower degree
neighbors of $t$ and $s.$ But $y$ or $z$ cannot be the common lower degree neighbor of $t$ and $s.$ Hence $d(y)=d(z)=2$ $\Rightarrow$ $d(s)\ge 3$ $\Rightarrow$
$s\leftrightarrow r$ and $ s\leftrightarrow t$ $\Rightarrow$ $d(r)=d(s)\ge 3.$ But now $r$ has no further neighbors besides $s,$ a contradiction.

The remaining case is that $t^{\prime}\nleftrightarrow z.$ Then $t^{\prime}\leftrightarrow x$ and $t^{\prime}\leftrightarrow y$ $\Rightarrow$ $d(x)=d(y)=2$ $\Rightarrow$ $z$ is the only possible
lower degree neighbor of both $t$ and $s$ $\Rightarrow$ $z\leftrightarrow t$ and $z\leftrightarrow s$ $\Rightarrow$ $d(z)=2$ and $d(s)\ge 3$ $\Rightarrow$
$s\leftrightarrow t$ and $s\leftrightarrow r$ $\Rightarrow$ $d(r)=d(s)\ge 3.$ But now $r$ has no further neighbors besides $s,$ a contradiction.

Subcase 2. Each neighbor of $v$ is not VT. Now $V(G)\setminus N[v]=\{t^{\prime},r,s\}$ where $t^{\prime},r$ and $s$ are VT, R and S respectively. Let $N(v)=\{a,b,c,t\}$ where $t$ is T.
We have $d(t)=3.$ First note that $t\nleftrightarrow t^{\prime}.$ Otherwise $d(t^{\prime})=2.$ But then $t^{\prime}$ has no lower degree neighbor. Hence $t$ can have
a lower degree neighbor only in $\{a,b,c\}.$ Next we distinguish three cases by considering where the third neighbor of $t$ with degree $3$ lies.

(1) $t\leftrightarrow s.$ Then $d(s)=3.$ We assert that $s\nleftrightarrow t^{\prime}.$ Otherwise $d(t^{\prime})=2.$ But then $t^{\prime}$ has no lower degree neighbor. Now
$N(t^{\prime})\subseteq \{a,b,c\}$ $\Rightarrow$ $2\le d(t^{\prime})\le 3.$ If $d(t^{\prime})=2, $ then $t^{\prime}$ has no lower degree neighbor, a contradiction;
if $d(t^{\prime})=3, $ then $t^{\prime}$ has no higher degree neighbor, a contradiction again.

(2) $t\leftrightarrow r.$ Then $d(r)=3.$ Since $r\nleftrightarrow t^{\prime},$ $r$ has at least one neighbor in $\{a,b,c\},$ say $r\leftrightarrow a.$ Consequently
$d(a)=3.$ Then $t$ can have a lower degree neighbor only in $\{b,c\},$ say $t\leftrightarrow b$ and $d(b)=2.$ Now $t^{\prime}$ can only have $c$ as its lower degree neighbor
$\Rightarrow$ $t^{\prime}\leftrightarrow c,$ $d(c)\ge 2$ and $d(t^{\prime})\ge 3$ $\Rightarrow$ $t^{\prime}\leftrightarrow a$ and $t^{\prime}\leftrightarrow s.$ Then $d(t^{\prime})=4,$ which is impossible since $\Delta =4$ and $t^{\prime}$ cannot have a higher degree neighbor.

(3) $N(t)\subseteq N[v].$ Without loss of generality, suppose $N(t)=\{v,b,c\},$ $d(b)=3$ and $d(c)=2.$ Then $t^{\prime}\nleftrightarrow b,$ since otherwise $d(t^{\prime})=4$
$\Rightarrow$ $t^{\prime}$ has no higher degree neighbor or $d(t^{\prime})=2$ $\Rightarrow$ $t^{\prime}$ has no lower degree neighbor . Now $N(t^{\prime})\subseteq \{a,s\}$ $\Rightarrow$ $d(t^{\prime})=2.$ But then $t^{\prime}$ has no lower degree neighbor,
a contradiction.  This completes the proof.$\Box$

\end{document}